\documentclass[11pt,reqno]{amsart}
\usepackage{txfonts}
\usepackage{a4wide}
\usepackage{amssymb}
\usepackage{caption}
\usepackage{epsfig}
\usepackage{subfigure}
\usepackage{graphics,color,psfrag}
\usepackage{float}

\allowdisplaybreaks[4]
\renewcommand{\baselinestretch}{1.2}

\newcommand{\al}{\alpha}
\newcommand{\be}{\beta}
\newcommand{\ga}{\gamma}
\newcommand{\eps}{\varepsilon}
\renewcommand{\th}{\vartheta}
\newcommand{\om}{\omega}
\newcommand{\percent}{\%}
\newcommand{\fC}{\mathbb{C}}
\newcommand{\fD}{\mathbb{D}}
\newcommand{\fH}{\mathbb{H}}
\newcommand{\fN}{\mathbb{N}}
\newcommand{\fQ}{\mathbb{Q}}
\newcommand{\fR}{\mathbb{R}}
\newcommand{\fT}{\mathbb{T}}
\newcommand{\fZ}{\mathbb{Z}}
\newcommand{\eg}{\textit{e.g.~}}
\newcommand{\ie}{\textit{i.e.~}}
\newcommand{\der}{\text{\textup{d}}}
\newcommand{\tle}{\preccurlyeq}
\newcommand{\tge}{\succurlyeq}
\newcommand{\tl}{\prec}
\newcommand{\tg}{\succ}
\newcommand{\p}{\partial}
\newcommand{\HH}{\mathcal{H}}
\newcommand{\DD}{\mathcal{D}}
\newcommand{\LL}{\mathcal{L}}
\newcommand{\RR}{\mathcal{R}}
\newcommand{\NN}{\mathcal{N}}
\newcommand{\VV}{\mathcal{V}}
\newcommand{\FF}{\mathcal{F}}
\newcommand{\GG}{\mathcal{G}}
\newcommand{\CC}{\mathcal{C}}
\newcommand{\TT}{\mathcal{T}}
\newcommand{\EE}{\mathcal{E}}
\renewcommand{\AA}{\mathcal{A}}
\newcommand{\BB}{\mathcal{B}}
\newcommand{\PP}{\mathcal{P}}
\renewcommand{\SS}{\mathcal{S}}
\newcommand{\VVVV}{\mathbf{V}}
\newcommand{\WWWW}{\mathbf{W}}
\newcommand{\defcol}[1]{\textit{#1}}
\newcommand{\io}{\infty}
\newcommand{\const}{\text{\textup{const.}}}
\newcommand{\epsc}{\varepsilon_{\text{c}}}
\newcommand{\changed}[1]{\textcolor[rgb]{0.8,0.2,0}{#1}}
\newcommand{\cutit}[1]{\textcolor[rgb]{0,0.4,0.4}{#1}}
\renewcommand{\captionfont}{\Small}
\newcommand{\bbe}{\betaup}

\renewcommand{\Im}{\operatorname{Im}}
\renewcommand{\Re}{\operatorname{Re}}
\newcommand{\am}{\operatorname{am}}
\newcommand{\sn}{\operatorname{sn}}
\newcommand{\cn}{\operatorname{cn}}
\newcommand{\dn}{\operatorname{dn}}
\newcommand{\Val}{\operatorname{Val}}
\newcommand{\sign}{\operatorname{sign}}
\newcommand{\Id}{\operatorname{\mathbf{1}}}

\theoremstyle{definition}

\theoremstyle{remark}

\newtheoremstyle{bgtheorem}{3pt}{3pt}{\itshape}{}{\bfseries}{}{ }%
{\thmname{#1}\thmnumber{ #2}\thmnote{ (#3)}}
\theoremstyle{bgtheorem}

\newtheoremstyle{case}{5pt}{5pt}{}{}{\scshape}{ }{ }%
{\thmnote{[#3]}}
\theoremstyle{case}

\numberwithin{equation}{section}

\begin{document}

\title[Scaling of the Critical Function]{
Scaling of the Critical Function for the Standard Map: \\
Some Numerical Results}

\author{Alberto Berretti}
\address{Alberto Berretti\\ Dipartimento di Matematica\\
II Universit\`{a} di Roma (Tor Vergata) \\
Via della Ricerca Scientifica, 00133 Roma, Italy\\
and Istituto Nazionale di Fisica Nucleare, Sez. Tor Vergata
}
\email{{\tt berretti@mat.uniroma2.it}}

\author{Guido Gentile}
\address{Guido Gentile\\ Dipartimento di Matematica\\ Universit\`{a}
di Roma Tre\\ Largo S. Leonardo Murialdo 1, 00146 Roma, Italy}
\email{{\tt gentile@mat.uniroma3.it}}

\begin{abstract}
The behavior of the critical function for the breakdown of the
homotopically non-trivial invariant (KAM) curves for the standard map,
as the rotation number tends to a rational number,
is investigated using a version of Greene's residue criterion.
The results are compared to the analogous ones
for the radius of convergence of the Lindstedt series,
in which case rigorous theorems have been proved.
The conjectured interpolation of the critical function
in terms of the Bryuno function is discussed.
\end{abstract}

\maketitle

\thispagestyle{empty}

\vspace{0.5cm}

\section{Introduction}\label{sect:intro}

A long-standing problem in the study of quasi-integrable Hamiltonian
systems is the characterization of the threshold for
the break-down of KAM invariant surfaces in terms of the arithmetic
properties of the frequencies vectors.  In this context, we consider a
simple, yet paradigmatic, discrete-time model, the so called
\emph{standard map}, introduced originally in \cite{Chi,Gr}.  The
standard map is the dynamical system defined by the iteration of the
map
\begin{equation}
     T_{\eps}: \; \begin{cases}
     x' = x + y + \eps\sin x \, ,\\
     y' = y + \eps\sin x \, .
     \end{cases}
        \label{eq:sm}
\end{equation}
Here $(x,y) \in \fT\times\fR$; but of course the map $T_{\eps}$ could
be lifted to a map
$$
T^{*}_{\eps}: (\xi,\eta) \mapsto (\xi',\eta')
$$
on the plane $\fR^2$ given by the same formula as \eqref{eq:sm}
with $(\xi,\eta)$ replacing $(x,y)$. For some background
information, we refer the reader to the enormous literature
on the topic, and in particular to \cite{MacKay1} for a review.

Despite its apparent simplicity, there are only a few properties of
the standard map which can be considered really well understood to
full extent, especially from an analytical point of view.  For
instance the existence of KAM invariant curves, for values of the
parameter $\eps$ small enough and Diophantine rotation numbers, has
been proved a long time ago, 
but only recently the dependence of the radius of convergence 
on the rotation number has been obtained \cite{D,BG2} as
an interpolation formula in terms of the Bryuno function (see below). 
Also for the studying of the separatrix splitting,
only recently the original program by Lazutkin \cite{L}
has been completely achieved in a rigorous way \cite{Gel}.

In particular no rigorous analysis has been implemented for detecting
the critical value of $\eps$ at which the KAM invariant curve breaks
down, and only numerical results and heuristic theories exist on that
subject; see \cite{MacKay1,MacKay3,AKW}.

In \eqref{eq:sm} we can eliminate the $y$ variable by writing the
dynamics ``in Lagrangian form'' as a second order recursion:
\begin{equation}
      x_{n+1}-2x_{n}+x_{n-1} = \eps\sin x_{n} \, ,
        \label{eq:recursion}
\end{equation}
for all $n\in\fZ$.

For $\eps=0$, the circles $y=(\const)$ are invariant curves on which the
dynamics is given by rotation with angular velocity $\om = y/2\pi$; we
call $\om$ the \emph{rotation number}.  Without generality loss we can
choose $\om\in(0,1)$ as the invariant curves of the standard map are
invariant under translation of $2\pi$ in the $y$-direction.

As the perturbation is turned on, we face the classical KAM problem of
determining which invariant curves survive and up to which size of the
perturbative parameter $\eps$.  Such invariant curves are given
parametrically by the equation
\begin{equation*}
        \mathcal{C_{\eps,\om}}:\;
        \begin{cases}
                x = \al + u(\al,\eps,\om) \, , \\
                y = 2\pi\om + u(\al,\eps,\om) - u(\al-2\pi\om,\eps,\om) \, ,
        \end{cases}
        \label{eq:conj}
\end{equation*}
where in the $\al$ variable the dynamics on the curve
$\mathcal{C_{\eps,\om}}$ is given by rotations $\al_{n+1} = \al_{n} +
2\pi\om$ (which solve \eqref{eq:recursion} for $\eps=0$).  The
function $u(\al,\eps,\om)$ is called the \emph{conjugating function}
or \emph{linearization}, and satisfies the functional equation
\begin{equation}
\begin{split}
      \left( D^{2}_{\om} u \right) (\al,\eps,\om)
      & \equiv u(\al+2\pi\om,\eps,\om) - 2u(\al,\eps,\om) +
      u(\al-2\pi\om,\eps,\om) \\
      & = \eps\sin(\al+u(\al,\eps,\om)) \, ,
      \label{eq:funceq}
\end{split}
\end{equation}
whose solutions are formally unique if we impose that
$u(\al,\eps,\om)$ has zero average in the $\al$ variable.  Therefore
the study of the invariant curves $\mathcal{C_{\eps,\om}}$ and of
their smoothness properties may be reduced to the study of the
existence and smoothness of the solutions of the functional equation
\eqref{eq:funceq}.

The solutions of \eqref{eq:funceq} can be studied perturbatively by
formally expanding $u(\al,\eps,\om)$ in Taylor series in $\eps$ and in
Fourier series in $\al$; the resulting series is what is traditionally
called the \emph{Lindstedt series}:
\begin{equation}
      u(\al,\eps,\om) = \sum_{k=1}^{\infty} \eps^{k} \, u^{(k)}(\al,\om) =
      \sum_{k=1}^{\infty} \eps^{k} \sum_{|\nu| \leq k}
      e^{i \nu \al} \, {u}^{(k)}_{\nu}(\om) \, .
      \label{eq:lind}
\end{equation}
To characterize the breakdown of an invariant curve
$\mathcal{C_{\eps,\om}}$ we introduce the \emph{radius of convergence}
of the Lindstedt series
\begin{equation}
      \rho(\om) = \inf_{\al \in \fT}
      \left(
      \limsup_{k \rightarrow \infty}
        \left|u^{(k)}(\al,\om)\right|^{1/k}
      \right)^{-1} \, ,
      \label{eq:rho}
\end{equation}
the \emph{lower (analytic) critical function}
\begin{equation}
      \epsc(\om) = \sup\{\eps' \geq 0: \;
      \forall\,\eps''<\eps'\;\mathcal C_{\eps'',\om}\;
      \text{exists and is analytic}\},
      \label{eq:epscrit}
\end{equation}
and the \emph{upper (analytic) critical function}
\begin{equation}
      \tilde{\eps}_{\text{c}}(\om) = \inf\{\eps' \geq 0: \;
      \forall\,\eps''>\eps'\;\mathcal C_{\eps'',\om}\;
      \text{does not exists as an analytic curve}\}.
      \label{eq:epscrit1}
\end{equation}
In general one could define analogous functions for negative values
of $\eps$; for the standard map they would be anyhow identical
(by symmetry properties).

Clearly $\rho(\om) \leq \epsc(\om)$ (in the early papers on the
subject some confusion was often made between $\rho$ and $\epsc$).  
It is instead \emph{believed}
that, \emph{for the standard map}, $\epsc(\om) =
\tilde{\eps}_{\text{c}}(\om)$, 
so we can speak generically of \emph{one} critical function $\epsc(\om)$
without further qualification. 
Note that for similar maps with
more general perturbations numerical results \cite{W} suggest that the
two critical functions may be indeed different.  Note also that one could
define breakdown thresholds with the analyticity condition in
\eqref{eq:epscrit}, \eqref{eq:epscrit1} replaced by a weaker one (such
as $C^{\infty}$ or $C^{k}$); again those thresholds could, in
principle, be different from the analytical one, though for the
standard map it is \emph{believed} that no such difference exists, so
that the analytic category is the right one to investigate the
breakdown phenomenon. 

The radius of convergence of the series \eqref{eq:lind}
is zero -- and no KAM invariant curve exists -- when $\om$ is
rational.  When $\om$ satisfies an irrationality condition known as
the \emph{Bryuno condition} (see below), instead, it can be proved
that $\rho(\om)>0$ -- so that analytic invariant curves exist for
$\eps$ small -- and even precise upper and lower bounds on the
dependence of $\rho(\om)$ on $\om$ can be given, up to a bounded
function of $\om$ \cite{D,BG2}.  More precisely for any rotation
number $\om$ one can define the \emph{Bryuno function} $B(\om)$,
as the solution of the functional equation \cite{Y}
\begin{equation}
     \begin{cases}
        B(\om) = -\log\om + \om B(\om^{-1}) \quad
            \text{for $\om \in (0,1)$ and irrational},\\
        B(\om+1) = B(\om) \, .
     \end{cases}
     \label{eq:bryuno}
\end{equation}
By an easy fixed-point argument it can be proved that a solution to
\eqref{eq:bryuno} exists and is unique in $L^{p}(\fT)$ for each $p \geq 1$.

We shall call \emph{Bryuno number} a number $\om$
satisfying the \emph{Bryuno condition} $B(\om)<\infty$.
Then for any Bryuno number $\om$ one has
\begin{equation}
        \left| \log \rho(\om) + 2 B(\om) \right| < C_{0} \, ,
        \label{eq:interpol}
\end{equation}
for a universal constant $C_{0}$, that is for a constant $C_{0}$
independent of $\om$ (see \cite{D,BG2} for a proof);
in particular this implies that an invariant curve with
rotation number $\om$ exists if and only if $\om$ satisfies
the Bryuno condition. Equation \eqref{eq:interpol} and similar formulas
are referred to as ``Bryuno's interpolation formulas''. 

The claim is often made that a formula analogous to
\eqref{eq:interpol} should hold for $\epsc(\om)$:
for any $\om$ satisfying the Bryuno condition one should have
\begin{equation}
        \left| \log \epsc(\om) + \beta B(\om) \right| < C_{1} \, ,
        \label{eq:interpol2}
\end{equation}
for a universal constant $C_{1}$, with an exponent $\beta \le 2$;
it is conjectured that $\beta=1$ (see \eg \cite{MS,CL}). 

Equation \eqref{eq:interpol2} implies a
scaling law for the critical function $\epsc(\om)$ as $\om
\rightarrow p/q$ \emph{on suitable sequences of Bryuno numbers}.  
In fact, given \eqref{eq:interpol2},
there are sequences of Bryuno (even Diophantine) rotation numbers
chosen in such a way that along them the critical function tends to
zero \emph{in any arbitrarily fast way}.  For example, we can consider
the two sequences of Diophantine (even noble) rotation numbers
\begin{equation}
    \om_{k} = \frac{1}{k+\ga}, \qquad
    \tilde{\om}_{k} = \cfrac{1}{k + \cfrac{1}{2^{k^{2}} + \ga}} \, ,
    \label{eq:twoseq}
\end{equation}
where $\ga = (\sqrt{5}-1)/2=[1^{\infty}]$ is the golden mean; then
\eqref{eq:interpol2}, with $\be=1$, would imply that $\epsc(\om_{k}) =
O(1/k)$ while $\epsc(\tilde{\om}_{k}) = O(e^{-k}/k)$, that is
\emph{much faster} (see \cite{BG2}, p.  625-626).  It is therefore
essential to have a good control over the arithmetic properties of the
rotation numbers one considers when speaking of scaling properties of
the critical function $\epsc(\om)$.

The conjecture of Bryuno's interpolation was actually made for the
critical function $\epsc(\om)$ more than 10 years ago in
\cite{MS}; in that paper, Bryuno's interpolation is stated 
formally for the radius of convergence, but the numerical
calculations,  with Greene's method, compute the critical function.
The main motivation behind \cite{MS} was the comparison with the work
of Yoccoz in \cite{Y}, together with the claims of universality coming
from the renormalization description of the critical invariant circle.
In \cite{D,BG2} (see also \cite{BG1}) Bryuno's interpolation for the
radius  of convergence was indeed proved; the mechanism of
proof in \cite{BG2}, based on the multiscale decomposition of the
propagators in the tree expansion, naturally generates an estimate of
$\rho(\om)$ in terms of the Bryuno function for the semi-standard and
standard maps.  On the other hand, there is no compelling \emph{a
priori} heuristic reason for the critical function $\epsc(\om)$ for
the standard map to satisfy an interpolation formula in terms of the
same arithmetical function as the radius of convergence $\rho(\om)$;
put it in another way, it is by no means obvious that  $|\log\rho(\om)
- (2/\be)\log\epsc(\om)|$ should be bounded. 

From this point of view, it would be interesting to consider
generalized standard maps, \ie maps where the nonlinear term in
\eqref{eq:sm} is an arbitrary analytic, periodic function of $x$ (see
\cite{BG3,BG4}). In these cases a Bryuno's interpolation formula for
the radius of convergence of their Lindstedt series is not known.

The method used in \cite{BFG} cannot be pushed so far to get reasonable
numerical data on the critical function, for some rather obvious
reasons; in fact, the method we used there (Pad\'{e} approximants)
attempts at modeling the \emph{whole} natural boundary, giving
particular weight at those regions of the boundary where the singularity
is ``strongest'': that is, to those regions closer to the origin (the
\emph{first order} or \emph{dominant singularities} as
defined in \cite{BFG}), which determine
$\rho(\om)$; so that part of the natural boundary near the real $\eps$ axis,
which determines $\epsc(\om)$, is represented, as $\om$ is closer and closer
to a rational value, as a few scattered points from which no reliable
information can be extracted: this happens already for rotation numbers
as little close to a rational value as, for instance, $1/(50+\ga)$ is
to $0$, that is still quite far from the rational value. One clearly
needs a method in which all the computing power is dedicated to the
calculation of the quantity one is interested in, that is $\epsc(\om)$.  

To this aim, two methods have been used previously in the literature:
Greene's method (also known as residue criterion; more about it in
the next section), used in \cite{MS}, and the frequency map analysis
\cite{La}, used in \cite{CL}. As we also use Greene's method, we
shall postpone a more thorough analysis to the next section, and go on
to a discussion of the results of \cite{CL}.

In \cite{CL} the following functions are defined:
\begin{equation}
    \om^{+}_{p/q}(\eps) = \inf\left\{\om > \frac{p}{q}: \CC_{\eps,\om}
    \text{ exists and is } C^{1}\right\},
    \label{eq:critompiu}
\end{equation}
and
\begin{equation}
    \om^{-}_{p/q}(\eps) = \sup\left\{\om < \frac{p}{q}: \CC_{\eps,\om}
    \text{ exists and is } C^{1}\right\}.
    \label{eq:critommeno}
\end{equation}
The meaning of those functions is that, for the given value of $\eps$,
no ($C^{1}$) invariant curves exist with rotation numbers between
$\om^{-}_{p/q}(\eps)$ and $\om^{+}_{p/q}(\eps)$.
The frequency map analysis method computes $\Delta_{p/q}(\eps) =
\om^{+}_{p/q}(\eps) - \om^{-}_{p/q}(\eps)$ for selected values of
$\eps$; $\Delta_{p/q}(\eps)$ should tend to $0$ with $\eps$ and in 
this way a lower bound on $\beta$ should be obtained (see below). Note
that $\eps$ is fixed, and correspondingly some rotation numbers are 
\emph{computed numerically}, therefore losing any strict control
over their arithmetical properties. 

We remark that the regularity properties of the functions
$\om^{\pm}_{p/q}(\eps)$ are quite hard to understand, and in particular
their relation with the critical function is far from obvious.
In fact, while it is certainly true that
\begin{equation}
    \epsc(\tilde{\om}) < \eps \quad \forall \tilde{\om} \in 
    (\om^{-}_{p/q}(\eps),\om^{+}_{p/q}(\eps)),
    \label{eq:fma1}
\end{equation}
the formulas at p. 2037 and p. 2052 of \cite{CL}, that is 
$\epsc(\om^{\pm}_{p/q}(\eps)) = \eps$, cannot be claimed in full rigor
since an invariant 
curve with rotation number very close to $p/q$ can be broken by the effect 
of \emph{another} resonance $p'/q' \approx p/q$, but distinct, so that
we can at most claim that
\begin{equation}
    \epsc(\om^{\pm}_{p/q}(\eps)) \le \eps.
    \label{eq:fma2}
\end{equation}
This implies that the law $|\om^{\pm}_{p/q}(\eps) -p/q| \approx
\eps^q$, numerically determined in \cite{CL}, provides for
the critical exponent an estimate from below of the actual value $\beta$,
which in principle could be higher (if it does exist at all). 
Equality in \eqref{eq:fma2} can be safely assumed at best for $\eps$ such that
the corresponding value $\om^{\pm}_{p/q}(\eps)$ belongs to
a special class of rotation numbers tending to $p/q$
(in some sense the ``best ones'', that is the ones whose
partial quotients grow as slow as possible), which are indeed
the ones considered in \cite{CL} and in the present paper
(and which are the only really accessible to a numerical investigation).
Note also that to saturate \eqref{eq:fma2} one should assume other 
qualitative features (like monotonicity) on the functions 
$\om^{\pm}_{p/q}(\eps)$, which are far from being proved. 
However \emph{for the case $\om \rightarrow 0$ only} this is enough, 
since estimates in \cite{TZ} imply \emph{an upper bound} 
on the critical exponent, which closes the gap. 

The numerical lower bounds found in \cite{CL} for $\beta$ are
consistent with $\beta=1$ with errors of orders
$4\percent$ for $\om$ close to $0/1$,
$10\percent$ for $\om$ close to $1/2$,
$5\percent$ for $\om$ close to $1/3$,
$10\percent$ for $\om$ close to $1/4$,
$8\percent$ for $\om$ close to $1/5$ and
$10\percent$ for $\om$ close to $2/5$.

Establishing a condition like \eqref{eq:interpol2}  is out of reach
from the numerical point of view if one wants to take into account
\emph{arbitrary} sequences of Bryuno numbers.  In fact for the
frequency map analysis this is a limitation intrinsic to the method
itself, since it automatically sort of chooses the best sequence of
Diophantine numbers tending to any given rational value. For any other
method, like Greene's residue criterion,  to investigate Bryuno
non-Diophantine numbers would require computer resources far beyond
current availability, while computing the critical function even for
Diophantine numbers with large partial quotients becomes substantially 
hard. So the question of establishing a Bryuno interpolation formula
for $\epsc(\om)$, and obtaining the correct critical exponent $\beta$
if such a formula is indeed established, is still quite open. 

In this paper we use Greene's method to compute the critical function
when the distance of the rotation numbers from the resonances is of
order $10^{-5}$. As the computations close to resonances become very
time-consuming we look at only three resonances ($0$, $1/2$ and $1/3$).
We then use the conjectured H\"older-continuity property of the
function $\log\epsc(\om)+\beta B(\om)$ to derive the corrections to the
asymptotic behavior of $\log\epsc(\om)$, so improving significantly
the agreement of the data with the conjectured value of $\beta=1$.

Of course the problem is not completely solved, even
from the numerical point of view, for two reasons.
The first is that we consider only three resonances,
so that a more exhaustive investigation would be needed.
The second is the aforementioned very special choice
of the sequences of rotation numbers tending to the resonances
that we have to use. Nevertheless we improve the results
existing in literature by one order of magnitude
both in the distance from the resonance and in the value
of $\beta$, finding further support for the
conjectured Bryuno's interpolation formula for the critical function.

\section{Greene's method}\label{sect:greene}

The main tool we use to determine numerically the break-down thresholds
for analytic invariant curves for the standard map is Greene's method,
known also as \emph{residue criterion}.  We now recall the main properties 
of the periodic solutions of the standard map used to formulate
Greene's method, and sketch briefly its foundations,
referring to the original paper \cite{Gr} for more details.

We also recall that in \cite{FL1} and \cite{MacKay2}
some theorems are proved that go some way in the direction of
proving the validity of Greene's method, at least in special cases.
While a full rigorous justification of its use has not yet been
achieved, Greene's method is considered one of the most accurate way
to compute the critical function $\epsc$ for the standard map. 

If $\om$ is a rational number, given as the irreducible fraction $p/q$,
then Birkhoff theory \cite{Bir2} applies; its consequences for maps
like the standard map $T_{\eps}$ are the following.  If $\eps = 0$
(unperturbed, linear case) then there are trivially invariant curves
with rational rotation number $p/q$, such that every point on them is
a fixed point of the iterated map $T^{\circ q}_{\eps}$.  As the
perturbation is turned on, only $2kq$, $k \in \fN$, points survive as
fixed points of the $q$-th iterate of the map $T_{\eps}$.  These
correspond to an even number ($2k$) of periodic orbits of period $q$.
Such orbits -- that we call \emph{perturbative} -- 
are the ones which will be studied within a perturbative
framework; a simple perturbative calculation (see \eg \cite{BG5})
shows that for the standard map the even number of such periodic
orbits is indeed just $2$.

Of course this does not mean that such orbits are the only periodic
ones for the standard map, but they are those which are obtained
by continuation (in $\eps$) from 
unperturbed ones.  
In other words such a scenario does not consider the new
periodic orbits arising when the perturbation is switched on.  If we
pass to the plane $\fR^{2}$ and consider the map $T^{*}_{\eps}$, then
the situation can be clarified in the following way.  When $\om$ is
irrational and satisfies the Bryuno condition, then the invariant
curve with rotation number $\om$ of the unperturbed map survives for
small values of $\eps$, while an invariant curve with rational
rotation number $p/q$ is suddenly destroyed; instead, only two
discrete invariant sets of points
$\{(\xi^{(\ell)}_{j},\eta^{(\ell)}_{j})\}_{j\in\fZ}$, $\ell = 1, 2$ 
survive, such that
\begin{equation}
\begin{cases}
    \xi^{(\ell)}_{j+1} > \xi^{(\ell)}_{j} \, , \\
    \xi^{(\ell)}_{j+q} = \xi^{(\ell)}_{j} + 2\pi p, \quad \ell = 1, 2.
\end{cases}
\label{eq:discreteinv}
\end{equation}
By taking the quotient in the first variable by the group of discrete
translations by multiples of $2\pi$, we obviously get two
periodic orbits of period $q$, on which the motion has rotation
number $p/q$.

In \cite{BG5} it is also proved that, for
small values of $\eps$, each such periodic orbit lies on an analytic
curve -- called a \emph{remnant} of the rational invariant curve of
the unperturbed map --, and for rational numbers which approximate a
Bryuno number $\om$ such remnants approximate the
invariant curve with rotation number $\om$.

The basic idea of Greene's method consists in relating the break-down of an invariant
curve with the loss of stability of nearby perturbative periodic orbits. 
In practice, the hypothesis behind Greene's method is that, if 
$\eps < \epsc(\om)$, then there is a sequence of stable perturbative 
periodic orbits with rotation numbers $p_{k}/q_{k}$; as $\eps$ grows 
beyond $\epsc(\om)$, these periodic orbits lose stability in the large
$k$ limit. 

The criterion can be formulated more precisely in the following way.
Let $\{(x^{(k)}_{i},y^{(k)}_{i})\}_{i=1}^{q}$ be a
perturbative periodic orbit with rotation number $p_{k}/q_{k}$,
approximating the irrational rotation number $\om$.
Let $\TT_{k}(\eps)$ be the trace of the tangent
dynamics along the periodic orbit:
\begin{equation}
        \TT_{k}(\eps) = \text{tr} \prod_{i=1}^{q_{k}}
        \begin{bmatrix}1+\eps\cos x^{(k)}_{i} & 1 \\
                \eps \cos x^{(k)}_{i} & 1
        \end{bmatrix}.
        \label{eq:trace}
\end{equation}
Then the periodic orbit is stable if $-2 < \TT_{k}(\eps) < 2$, unstable
otherwise.  For historical reasons, the criterion is usually formulated
in term of the \emph{residue} $\RR_{k}(\eps)$ of the orbit, 
related to the above trace by
$$
\RR_{k}(\eps) = \frac{2-\TT_{k}(\eps)}{4}.
$$
Therefore in terms of the residue the orbit is stable if
$0 < \RR_{k}(\eps) < 1$, unstable otherwise.
We then track, for a fixed value of $\eps$, the residue of
those perturbative periodic orbits with rotation numbers $p_{k}/q_{k}$
which are stable for $\eps=0$; if the residue diverges as
$k\to\infty$, then $\eps > \epsc(\om)$,
while if the residue tends to $0$ then $\eps < \epsc(\om)$. 

It is actually conjectured (see \cite{MacKay2}) that
if $\eps < \epsc(\om)$, then the residue $\RR_{k}(\eps)$
tends exponentially to zero as $k \rightarrow \infty$,
with a rate of decay proportional to the width of the analyticity 
strip of the conjugating function $u(\al,\eps,\om)$ on the complex $\al$ 
plane for the values of $\eps$ and $\om$ considered. So Greene's method 
can also be used also to provide numerical information on the analytic 
properties of $u$ in $\al$, assuming this conjecture. 

An interesting question is what actually happens to the residue
\emph{at the critical function} $\epsc(\om)$.  It was originally
conjectured that for \emph{noble} rotation numbers
(that is, rotation numbers which are obtained by
applying a modular transformation to the golden mean, so that their
continued fraction expansion has a ``tail'' of $1$'s)
it tends to a limit value, which should be about $0.25$.
We present below some numerical results which show that
generally the situation is more complicate,
and that such \emph{limit residue} $\RR_{\infty}(\eps)$
could be not only different for different classes of rotation numbers,
but could also be non-existent, and relate the behavior of the
sequence of residues $\RR_{k}(\eps)$ for a fixed value of $\eps$ along
the sequence of perturbative periodic orbits of rotation numbers
$p_{k}/q_{k}$ to the arithmetic properties of the rotation number. 

From the practical, computational standpoint, the implementation of
Greene's method faces some challenges if we wish to use it near
resonances.  The first is that, if $\om$ is near a resonance, then the
$q_{k}$ become soon quite large, that is we have to find many
\emph{long} periodic orbits, which takes a lot of computer time.  

The second, hardest, challenge is more subtle.  In fact, if it happens
that the rotation number of a periodic orbit is $p/q \approx p'/q'$,
with $q \gg q'$ (the typical situation arising when approximating
irrational rotation numbers close to small-denominator rationals) then
it appears numerically that the periodic orbit of rotation number $p/q$
tends to consist in lots ($q$ is supposed to be large) of points
accumulating near the points making the periodic orbit of rotation
number $p'/q'$.  The consequences for the computation of the residue
are dire, as in this case the matrix in \eqref{eq:trace} has two very
large, opposite, nearly equal in absolute value diagonal elements, so
that when computing the trace the real data cancels and one is left
with just the numerical error. Note that using a low precision
with Greene's method so close to resonances gives essentially
noise instead of the residue, so we get no values at all for $\beta$.
We choose a brute-force solution to this precision problem,
which consists in increasing the number of digits in the
calculations until some data is left when computing the trace.
Empirically, this could mean that one has to use \emph{hundreds} of
digits of precision in computing \eqref{eq:trace} numerically,
\emph{therefore also the periodic orbits must be known with such a
precision}: considering that one easily needs periodic orbits of
period in the range of several tens of thousands -- we actually reach
orbits of length of the order of $150000$ --, the calculation of a
single value of $\epsc$ can require a great amount of computer time.

\section{Numerical results}\label{sect:results}

\subsection{Rotation numbers close to 0}\label{subsect:zero/one}

Consider rotation numbers $\om_{n}=1/(n+\ga)=[n,1^{\infty}]$,
with $n\in \fN$: in table \ref{tab:1} we give the values of
the Bryuno function and of the critical function for rotation numbers
$\om_{n_{k}}$, with $\{n_{k}\}$ a finite increasing sequence.
Note that we reach values of rotation number close
more than $2\times 10^{-5}$ to the resonance value ($0$ in this case),
which corresponds to values of $n$ up to $60000$.

By fitting $y=-\log\epsc(\om_{n_{k}})$ as a linear function of
$x=B(\om_{n_{k}})$, we obtain
\begin{equation}
        y= ax + b \qquad a = 0.9705 , \qquad b = -1.9553 .
        \label{eq:bestfit1}
\end{equation}
As we see the slope is close to (but different from) 1:
the relative difference is about $3.0\percent$.

One also realizes that the slope of the line increases if
we neglect the rotation numbers $\omega_{n}$ corresponding 
to smaller values of $n$: this suggests that, if we consider just
pairs of successive rotation numbers and evaluate the slope
of the line passing through them, then we obtain an increasing function.
This can be formulated more precisely as follows.
For $n,m\in \fN$ define
\begin{equation}
        A(\om_{n},\om_{m}) = - \frac{\log \epsc(\om_{n})
	- \log \epsc(\om_{m})} {B(\om_{n})-B(\om_{m})} ,
        \label{eq:slope}
\end{equation}
which measures the slope $a$ of the line
\begin{equation}
        -\log\epsc(\om)= a B(\om) + b ,
        \label{eq:line}
\end{equation}
passing through the points $(B(\om_{n}),-\log
\epsc(\om_{n}))$ and $(B(\om_{m}),-\log \epsc(\om_{m}))$.
We set $A_{k}=A(\om_{n_{k+1}},\om_{n_{k}})$.
In table \ref{tab:1} we give also the values of the slopes $A_{k}$:
as we noted $A_{k}$ steadily increases.

The value $\beta=1$ is anyhow still far from being reached: at best,
just considering the last value of $A_{k}$ in table \ref{tab:1} we
obtain a value whose relative difference from 1 is greater than
$1\percent$. Moreover, though the values of the slopes increase
as $n\to\io$, the convergence to 1 is very slow.  In the next section
we shall provide a heuristic argument which allows  to guess the
correction to the asymptotic behavior and so try to extrapolate a
better value of $\beta$; this applies also to the cases considered in
the next subsections. 

In table \ref{tab:2} we give the values of the Bryuno function
and of the critical function for a finite sequence
of rotation numbers $\om_{n_{k}}$, with
$\om_{n_{k}}=1/(n_{k}+1/(20+ \gamma))=[n_{k},20,1^{\infty}]$:
such numbers tend to $0$ as the previously considered ones,
and share with them, essentially, the same Diophantine properties,
as they have the same ``tail'' of $1$'s in their continued fraction
expansion, with the only difference that
there is a partial quotient $20$ before such a ``tail''. 
The distance of the rotation numbers considered is up to
$2\times 10^{-4}$ from $0$, i.e. an order less than in the previous
case: this is due to the fact the partial quotients
go faster, and it becomes longer for the residue to reach
the asymptotic value (so that periodic orbits
with larger periods should be considered in order to obtain
for the rotation numbers the same distance from the resonance value).

As one can see, the values of the Bryuno function and of
the critical function are comparable with those listed
in table \ref{tab:1}: the introduction of a larger partial quotient 
does not introduce any relevant change. As a consequence,
also the slopes $A_{k}$, defined as before with the new
definition of $\om_{n_{k}}$, are very similar 
(as a look at the last column of table \ref{tab:2} immediately confirms).

Note however that to compute numerically the critical function
for rotation numbers of the form $[n_{k},20,1^{\infty}]$
for given $k$ is much more time consuming, since, in general,
to obtain a reliable precision we are forced to reach periodic orbits
with very high periods (say more than a hundred thousand),
which requires a precision of about 600 digits.

\subsection{Rotation numbers close to 1/2}\label{subsect:one/two}

In table \ref{tab:3} we consider a sequence of rotation numbers tending
to $1/2$ of the form $\om_{n}=1/(2+1/(n+\gamma))=[2,n,1^{\infty}]$.
The rotation numbers considered are up to $10^{-5}$ close
to the resonance value $1/2$ (which correspond to values of $n$ up
to 20000).

The fit for $y=-\log\epsc(\om_{n_{k}})$ as a linear function of 
$x=B(\om_{n_{k}})$ gives
\begin{equation}
        y= ax + b \qquad a = 0.9641 , \qquad b = -1.6203 .
        \label{eq:bestfit2}
\end{equation}

Again we see the slope is not 1, and
the relative error is now about $3.6\percent$.
It is greater than in the previous case because we stopped to
smaller values of $n$; in fact the values of the slopes
listed in table \ref{tab:3} show that again the
function $A_{k}$, defined exactly as before with the new
definition for the rotation numbers $\om_{n_{k}}$, is increasing in $k$.
The relative difference from 1 of the last value of $A_{k}$
is about $1.7\percent$.

\subsection{Rotation numbers close to 1/3}\label{subsect:one/three}

In table \ref{tab:4} we consider a sequence of rotation numbers tending
to $1/3$ of the form $\om_{n}=1/(3+1/(n+\gamma))=[3,n,1^{\infty}]$. The
rotation numbers considered are up to $5 \times 10^{-6}$ close to the
resonance value $1/3$ (which correspond to values of $n$ up to 20000).

The fit for $y=-\log\epsc(\om_{n_{k}})$ as a linear function of 
$x=B(\om_{n_{k}})$ gives
\begin{equation}
        y= ax + b \qquad a = 0.9637 , \qquad b = -1.6526 .
        \label{eq:bestfit3}
\end{equation}
while in the last column of table \ref{tab:4} we list the slopes
$A_{k}$, again defined as before with the new definition for the
rotation numbers $\om_{n_{k}}$; the relative difference of the slope
with respect to 1 is more than $3.7\percent$, while
the relative difference from 1 of the last value of $A_{k}$
is about $1.7\percent$.

\subsection{Behavior of the critical residues and other rotation numbers}

The behavior of the residue for $\eps$ \emph{exactly} equal to the
critical function $\epsc(\om)$ when the rotation number of the
approximating periodic orbits tends to $\om$ has been considered since
the very first papers on the subjects (for example, in \cite{Gr}
itself). In particular, one considers the sequence of residues
$\RR_{k}(\epsc(\om))$ when $k \rightarrow \infty$; it appears that this
sequence has a limit only when $\om$ is a number of so called ``constant
type'', \ie when $\om$ can be written as
$[a_{1},\dots,a_{N},d^{\infty}]$. This limit moreover seems to depend
only on the integer $d$, and not at all from the ``head'' of the
continued fraction expansion $[a_{1},\dots,a_{N}]$. Unfortunately, a
sound numerical evidence can be obtained only for $d=1$ and for short
``heads'' in the continued fraction expansion, otherwise the partial
quotients gets soon large and it becomes difficult to compute the
critical residues with the accuracy required: therefore we state this
more as a somewhat numerically founded and reasonable conjecture than 
else. In table \ref{tab:5} we give some values of the critical residue
for a few values of $d$. 

If a rational number is not of constant type, then a limit does not seem
to be achieved for the sequence of critical residues. In fact, it seems
to happen that if $\om$ is a quadratic irrational, so that the sequence
of the partial quotients $a_{k}$ is eventually periodic, the sequence of
critical residues is itself eventually periodic with the same period. In
tables \ref{tab:6}, \ref{tab:7}, \ref{tab:8} we can see the sequence
of critical residues for some quadratic irrationals with short periods
(resp. $2$, $2$ and $3$). If the rotation number is not a quadratic
irrational, so the partial quotients are aperiodic, the sequence of
critical residues does not seem to have any regularity (but see below
for a numerical difficulty). 

So far only quadratic irrational $\om$ have been considered. This is of
course a limitation, due mainly to practical reasons; in fact, quadratic
irrationals are the only irrationals with an eventually periodic
continued fraction expansion, so they are particularly suited to
Greene's residue criterion for two reasons: (1) the partial quotients
$a_{k}$ are bounded, since they are periodic, so the approximants
$p_{k}/q_{k}$ have denominators which do not grow too much and (2) if
the period is reasonably small, one can tell whether the critical
function has been reasonably approximated by looking at the sequence of
residues over a span of periods and easily see whether it decreases or
increases instead of being periodic. Instead, if the sequence of the
partial quotients is aperiodic (and worst yet, unbounded) one can never
be sure that the critical function has been obtained since the next
periodic orbit to be considered (corresponding to the next approximant
$p_{k}/q_{k}$) could come from an abnormally high (or low) partial
quotient $a_{k}$. Note that in \cite{CJ}, where general irrationals
are also considered in Subsection 3.3, in the numerical calculations
of the critical function, only the first ten partial quotients
of the rotation numbers are retained, and all the others are set to 1,
so one practically comes back to the case of noble numbers like ours.

This raises the question whether the \emph{algebraic} (rather
than just number-theoretic) properties of the
rotation number have any role in the properties of the corresponding
invariant curve. Lindstedt series expansion methods for instance do not
care about the algebraic properties of $\om$, as the only relevant
property is whether $\om$ is a Bryuno number or not. ``Phase space''
renormalization group methods instead seem to work (or at least they
have been applied) only in the case of quadratic irrational rotation
numbers, so their results could depend on the algebraic layer. We expect
that the algebraic properties of $\om$ could show up, maybe, in
discussing the smoothness of the natural boundary in the complex $\eps$
plane, but of course this is just a speculation.

\section{Discussion}\label{sect:discussion}

Despite intensive numerical calculations, the problem of confirming
the  conjecture expressed in \eqref{eq:interpol2} and estimating the
exponent $\beta$ cannot be considered completely settled even from the
numerical point of view. In fact, as we noted earlier, only the three
resonances $0/1$, $1/2$ and $1/3$ have been considered, and only very
special sequences of rotation numbers tending to such rational numbers
have been used: considering other sequences of rotation numbers, in
fact, means using numbers which have quite soon very large partial
quotients, so that they are very bad from the numerical point of view. 

Moreover, a simple linear fit of $\log\epsc(\om)$ against $B(\om)$,
that is a fit which takes into account only the leading conjectured
asymptotic behavior without any corrections, still gives results which
are quite unsatisfying, as the difference between the estimated value
of $\beta$ and the conjectured value $\beta=1$ is still of the order of
a few percent. What is worst, the ``running slopes'' $A_{k}$ defined in
the preceding section continue to grow monotonically from below, slowly
but steadily, so that one cannot even conclude that the conjecture is
false or that the value of $\beta$ is actually smaller than $1$. Clearly
corrections must be taken into account, or otherwise rotation numbers
even closer to the resonances (and significantly such) must be considered,
which is numerically unfeasible with current resources. 

Note also the apparently quite singular fact that for $\rho(\om)$ the
value $2$ of the  corresponding critical exponent seems to be obtained
within a few percent \emph{much earlier}. For instance for the rotation
numbers $\om$ close to $1/2$ listed in table \ref{tab:9}, by using for
the corresponding radii of convergence the values $\rho_{P}(\om)$
computed by Pad\'{e} approximants, we obtain for the slopes
$A_{k}'=A'(\om_{n_{k+1}},\om_{n_{k}})$, with
\begin{equation}
        A'(\om_{n},\om_{m}) = - \frac{\log \rho(\om_{n}) -
	\log \rho(\om_{m})} {B(\om_{n})-B(\om_{m})} ,
        \label{eq:rhoslope}
\end{equation}
the values in the last column of table \ref{tab:9}. Analogously for the
rotation numbers $\om$ close to $1/3$ listed in table \ref{tab:10},
again by using the values $\rho_{P}(\om)$ computed by Pad\'{e}
approximants for the corresponding radii of convergence, we obtain the
slopes $A_{k}'$ in the last column of the same table.

\begin{figure}[h]
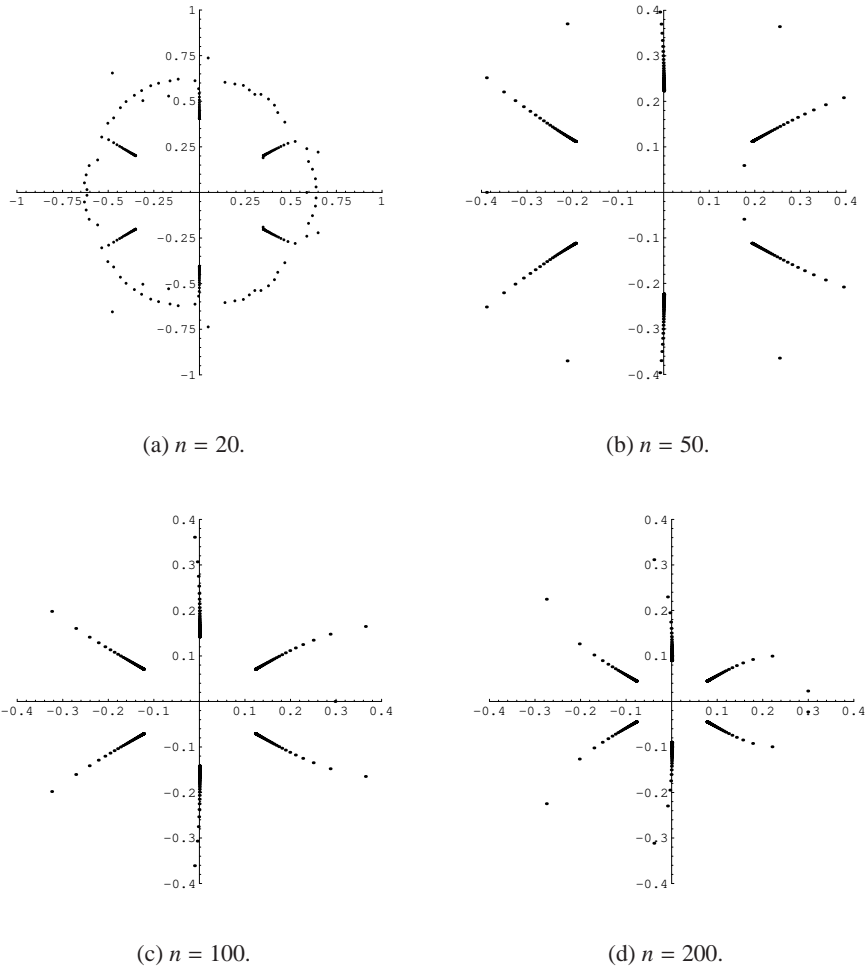

        \centering
  \subfigure[$n=20$.]{
    	\psfig{file=bg7fig1a.eps,height=2.0in}}
     \hspace{0.3 in}
  \subfigure[$n=50$.]{
    	\psfig{file=bg7fig1b.eps,height=2.0in}}
  \subfigure[$n=100$.]{
    	\psfig{file=bg7fig1c,height=2.0in}}
     \hspace{0.3 in}
  \subfigure[$n=200$.]{
    	\psfig{file=bg7fig1d,height=2.0in}}
        \caption{Poles of the Pad\'{e} approximant [240/240] for
	$\om=[3,n,1^{\infty}]$ and $\al=1$.
        }
        \label{fig:1}
\end{figure}

In figure \ref{fig:1} we represent the analyticity domains for
$\om=[3,20,1^{\infty}]$, $[3,50,1^{\infty}]$, $[3,100,1^{\infty}]$ and
$[3,200,1^{\infty}]$ as given by the poles of the Pad\'{e} approximants
$[240/240]$. As noted in Section \ref{sect:intro} for $\om$ getting
closer to $1/2$ the poles tend to accumulate near the strongest
singularity: therefore Pad\'{e} approximants are not suitable for
determining the critical function, but they can be fruitfully used in
order to detect the radius of convergence.

The relative errors with respect to $2$ for the values corresponding to
$n=40$, taken from tables \ref{tab:9} and \ref{tab:10},
are about $1.4\percent$ and $3.0\percent$, respectively,
therefore they are comparable with the errors for the last entries of
the corresponding tables \ref{tab:3} and \ref{tab:4} for the critical
function for rotation numbers \emph{much closer} to the resonance
values: in the latter case indeed such errors are about $1.7\percent$.
And for larger values of $n$
the relative errors become much smaller: for instance,
for $n=100$ and $n=200$, we find from table \ref{tab:10}
errors about $1.0\percent$ and $0.5\percent$, respectively.

Of course it would be also interesting to have the slopes for the
rotations numbers appearing in tables \ref{tab:3} and \ref{tab:4}. To
obtain the values of $\rho(\om)$ numerically can be as hard as to
determine the critical function $\epsc(\om)$. Also using the method of
Pad\'e approximants can be delicate,  as in order to obtain reliable
results a very high precision could be necessary. One could think of
using the complex extension of Greene's method envisaged in \cite{FL2},
and the analysis, at best, could be as delicate as in the present
paper, where real values of $\eps$ have been studied. We have also two
more difficulties with respect to the case of $\epsc(\om)$. First one
has to guess the direction in the complex plane where the singularities
of the boundary of the analyticity domain are the closest to the
origin; in this respect the results of \cite{BFG} suggest, as a natural
Ansatz, that, for rotation numbers close to $p/q$,
they such singularities lie along the directions of the
$2q$th roots of $-1$. Next, for fixed $\om$ close to a resonance value,
the value of $\rho(\om)$ should be much smaller than the value of
$\epsc(\om)$, again as a byproduct of the numerical analysis of
\cite{BFG} (and also that of \cite{CL}), so that the value of
$\rho(\om)$ is expected to be harder to detect than $\epsc(\om)$, as it
should require more precision and hence more computing time. However we
prefer to avoid any technical difficulties and to circumvent the
problem by using the heuristic formula introduced in \cite{BFG}, say
$\rho(\om) \approx \rho_{1}(\om)$, with
\begin{equation}
 	\rho_{1}(\om) = \eta^{2/q}
	\left( q |C_{p/q}|^{-1} \lambda_{c} \right)^{1/q},
	\label{eq:rhoapprox}
\end{equation}
where $\eta=|\om-p/q|$ if $\om$ is close to the resonance $p/q$,
$C_{p/q}$ is the numerical constant introduced in \cite{BG1} (one has
$C_{0/1}=1$, $C_{1/2}=-1/8$ and $C_{1/3}=-1/24$), and
$\lambda_{c}=4\pi^{2}\times 0.827524 \approx 32.669338$.

In \cite{BFG} we have already seen that there is a good agreement
between the value $\rho_{P}(\om)$ of the radius of convergence found by
Pad\'{e} approximants and the value $\rho_{1}(\om)$ predicted by the
formula \eqref{eq:rhoapprox}. Furthermore the formula
\eqref{eq:rhoapprox} becomes more and more reliable as $\om$ approaches
an rational value. See for instance tables \ref{tab:9} and
\ref{tab:10}, which show how the difference between the two values
$\rho_{P}(\om)$ and $\rho_{1}(\om)$ tend to shrink to zero when making
the rotation number $\om$ closer to the rational values $1/2$ and
$1/3$, respectively. So we can expect that the approximation we make by
evaluating the radius of convergence $\rho(\om)$ with $\rho_{1}(\om)$
is very good for values much closer to the resonance values, as the
ones we have considered are.

Then we obtain the values listed in tables \ref{tab:11}, \ref{tab:12}
and \ref{tab:13} for values of $\om$ close, respectively, to 0, 1/2 and
1/3 (the same for which we determined numerically the critical
function); the slopes $A_{k}'$ are listed in the last columns of these
tables.  Of course, if we use the formula \eqref{eq:rhoapprox}, a slope
approximately equal to $2$ is expected, by the definition itself of
$\rho_{1}(\om)$. The important fact is, in any case, that the
discrepancy with respect to the value $\beta=2$ (which in such a case
\emph{is known} to be the right one) is much smaller. In other words
the asymptotic formula \eqref{eq:interpol} is reached much earlier than
the one which is believed to hold for the critical function.

This different speed in reaching the asymptotic behavior of
$\rho(\om)$ and $\epsc(\om)$  can be explained in terms of different
corrections to the leading order when $\om \rightarrow p/q$ (and
therefore $B(\om) \rightarrow \io$). We shall now try to compute such
correction, at least heuristically, both for $\rho(\om)$ and for
$\epsc(\om)$, and try to use them to extrapolate a better value of
$\beta$. 

For what concerns $\rho(\om)$, we shall assume the validity of the 
heuristic formula \eqref{eq:rhoapprox}; this of course can introduce 
further corrections not accounted for, which we neglect assuming that 
they are smaller. 

Consider for example $\om_{n} = 1/(n+\ga)$, $n \rightarrow \io$. Then 
$\eta = \om_{n}$ in \eqref{eq:rhoapprox} and
\begin{equation*}
B(\om_{n}) = -\log\om_{n} + \om_{n}B(1/\om_{n}),
\end{equation*}
which implies that $\log(n+\ga) = B(\om_{n}) - B(\ga)/(n+\ga)$. 
Therefore in first approximation we have that
\begin{equation*}
\log(n+\ga) \approx B(\om_{n}) - B(\ga)e^{-B(\om_{n})}, 
\end{equation*}
as the leading behavior of $B(\om_{n})$ is just 
$\log(n+\ga)$. This gives
\begin{equation}
 	\log\rho_{1}(\om_{n}) \approx 
	\log(|C_{0/1}|^{-1} \lambda_{c})
	-2B(\om_{n}) + 2B(\ga) e^{-B(\om_{n})},
	\label{eq:rhocorr1}
\end{equation}
that is the correction to the linearly growing asymptotic behavior is
\emph{exponentially small}. An analogous, slightly more complicated, 
computation for $\om_{n} = 1/(q+1/(n+\ga))$ gives a correction of the
form  $B(\om_{n})\exp(-qB(\om_{n}))$, that is still essentially
exponentially small.  This explains the exceptional rapidity of the
approach to the scaling behavior  for $\rho(\om)$. 

Quantitatively, a fit of the numerical data of table \ref{tab:11} using
\eqref{eq:rhocorr1} to model the data gives for $\beta$ the value
$1.9999989$, whose difference from the correct value of $2$ is of the
order of $10^{-7}$, while a straight linear fit gives $2.00091$, whose
error is three orders of magnitude larger. 

Analogously, a fit of the numerical data of table \ref{tab:13} using
\begin{equation}
 	\log\rho_{1}(\om_{n}) \approx 
	\log(|C_{0/1}|^{-1} \lambda_{c})
	- \beta B(\om_{n}) + \left(b + c B(\om) \right) ^{-3B(\om_{n})},
	\label{eq:rhocorr3}
\end{equation}
gives for $\beta$ the value $2.000000287$,
whose difference from the correct value of $2$ is of the
order of $10^{-7}$, while a straight linear fit gives
$1.99984$, whose error is again three orders of magnitude larger.
Also a comparison between the mean-square distances of the data
from the corresponding fits is remarkable: we obtain
$2.389\times 10^{-8}$ for the fit by \eqref{eq:rhocorr3},
and $0.0000858$ for the linear fit. Note that also fits with
either $b=0$ or $c=0$ in \eqref{eq:rhocorr3} are worse:
for $b=0$ we obtain a value $\beta=1.99966$
(with mean-square distance $0.0000413$), while for $c=0$
we obtain a value $\beta=1.99967$ with mean-square distance $0.0000393$).

To compute the correction to the leading behavior of $\epsc(\om)$ is
of course quite another matter, since we don't even have a proof or at
least a very strong theoretical argument for the leading order. So the
following argument is more a qualitative explanation rather than a
quantitative attempt to extrapolate seriously the value of $\beta$ (but
we shall try nevertheless). 

As before we shall consider only the case 
$\om_{n} = 1/(q+1/(n+\ga))$, $n \rightarrow \io$, and we shall set
$\eta_{n}=|\om_{n} -1/q|$. Let
\begin{equation}
	\log\epsc(\om_{n}) + B(\om_{n}) = C(\om_{n}),
	\label{eq:epscorr1}
\end{equation}
where the function $C(\om)$ is believed to be
continuous (see \cite{Marmi,MS}). Let then  $\bar{c}
= \lim_{\om\rightarrow 0} C(\om)$.
It is also conjectured (see \cite{MMY}) that $C(\om)$ is
H\"older-continuous with some exponent $\al$
(in the quoted paper it is suggested that $\al$ could be $1/2$):
so, by recalling that $\eta_{n} \approx e^{-q B(\om_{n})}
\rightarrow 0$, a reasonable guess in \eqref{eq:epscorr1} could be
\begin{equation}
	\log\epsc(\om_{n}) = (\const) - \beta B(\om_{n}) +
	O(e^{-\al q B(\om_{n})}) .
	\label{eq:epscorr5}
\end{equation}
Despite the rough, qualitative nature of the argument above, we can try
to fit the data with the formula \eqref{eq:epscorr5} and see whether
the  growth of the slopes is such that the value of $1$ can actually be
reached. The ``best'' value of $\alpha$ is obtained by choosing it in
such a way that the mean square distance of the experimental data from
the values obtained from the fit is minimal. As an alternative, we
performed also nonlinear fits using Levenberg-Marquardt method (see
\cite{Leven,Marq}), obtaining consistent results. 

\begin{figure}[h]
\centering
  \subfigure[Sequence $\om_{n_{k}}$ listed in table \ref{tab:1}.]{
    	\psfig{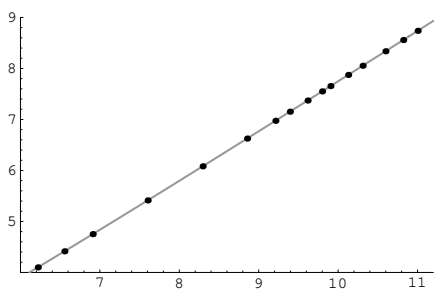}}
\hspace{0.3in}
  \subfigure[Sequence $\om_{n_{k}}$ listed in table \ref{tab:3}.]{
    	\psfig{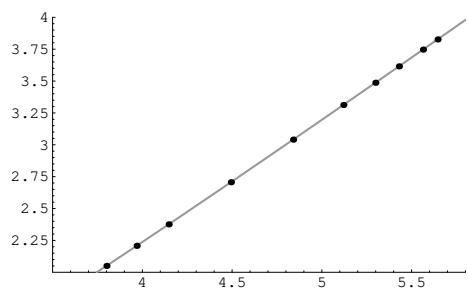}}
  \subfigure[Sequence $\om_{n_{k}}$ listed in table \ref{tab:4}.]{
    	\psfig{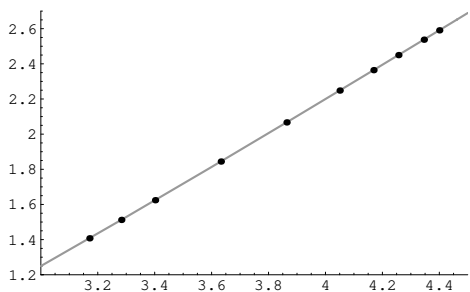}}
        \caption{Numerical values of $-\log\epsc(\om_{n_{k}})$,
	obtained with Greene's method, versus $B(\om_{n_{k}})$
        for the sequence $\om_{n_{k}}$ listed in tables \ref{tab:1},
	\ref{tab:3}; and \ref{tab:4};
	the error bars are less than the size of the points.
	The solid curve corresponds to the fits \eqref{eq:epsfit1},
	\eqref{eq:epsfit2} and \eqref{eq:epsfit3}.
        }
        \label{fig:2}
\end{figure}

Fitting the data relative to $\epsc(\om_{n})$ for the sequence
considered in table \ref{tab:1} with the formula \eqref{eq:epscorr5},
we obtain
\begin{equation}
	\log\epsc(\om_{n}) \approx
	-2.34630 +  1.00359 \, B(\om_{n}) +
	1.59684 \, e^{-0.3302 \, B(\om_{n})} ,
	\label{eq:epsfit1}
\end{equation}
that is finally a value much closer to $1$ than the straight linear fit,
which gave  $0.97052$. Moreover, the mean-square distance of the data
from the fit is $0.000210$ in the case of the fit with corrections,
while is much larger, that is $0.0396$, in the case of the linear
fit (see figure \ref{fig:2}a). 

If we consider the case $\om \rightarrow 1/2$ as in subsection
\ref{subsect:one/two}, the fit with formula \eqref{eq:epscorr5} gives
\begin{equation}
\log\epsc(\om_{n}) \approx 
	-1.86364 + 1.00308\,B(\om_{n}) +
	1.43766 \, e^{-0.69671\,B(\om_{n}) } 
\label{eq:epsfit2}
\end{equation}
with mean-square distance $d=0.0000512$
(the linear fit would give $\beta=0.96413$ and $d=0.0124$);
in figure \ref{fig:2}b we plot the data together with the fit. 

Finally, in the case $\om \rightarrow 1/3$ as in subsection
\ref{subsect:one/three}, the fit with formula \eqref{eq:epscorr5} gives
\begin{equation}
\log\epsc(\om_{n}) \approx 
                  -1.84393 + 1.00344\,B(\om_{n}) +
		  1.82643 \, e^{-1.0300 \, B(\om_{n}) },
\label{eq:epsfit3}
\end{equation}
with mean-square distance $0.0000403$
(the linear fit would give $\beta=0.96369$ and $d=0.00832$);
in figure \ref{fig:2}c we plot the data together with the fit. 

The results, together with the interpolation formula
\eqref{eq:epscorr5}, hint at a value of $\al$ close to $1/3$,
while $C_{\rho}(\om)$, according to the formula \eqref{eq:rhocorr3},
seems to be H\"older-continuous with any
exponent $\al_{\rho}<1$ (and $\al_{\rho}=1$ in $0$). So while in both
cases the corrections are exponentially small in $B(\om)$, the
coefficient in the exponential is about three times larger for
$\log\rho(\om)$, leading to smaller corrections and faster approach to
the asymptotic regime. 

\begin{figure}[h]
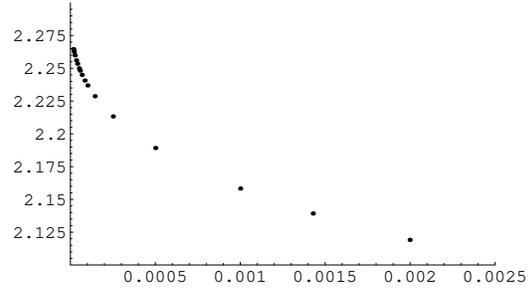
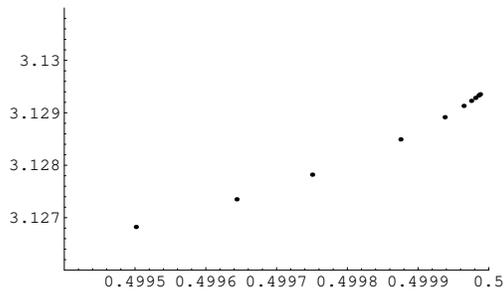
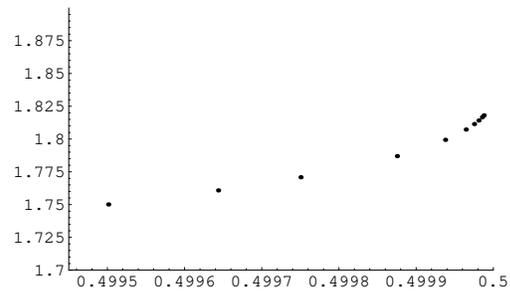
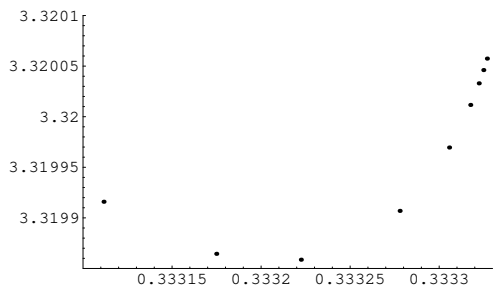
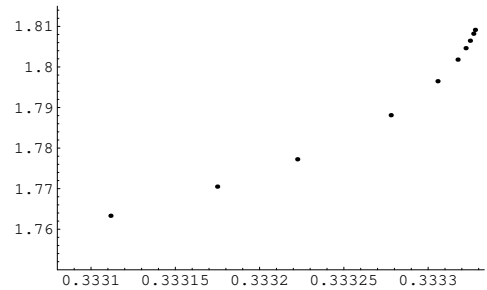

        \centering
  \subfigure[$C_{\rho}(\om)$ for $\om$ close to $0$.]{
    	\psfig{file=bg7fig3a.eps,height=1.5in}}
     \hspace{0.3 in}
  \subfigure[$C(\om)$ for $\om$ close to $0$.]{
    	\psfig{file=bg7fig3b.eps,height=1.5in}}
  \subfigure[$C_{\rho}(\om)$ for $\om$ close to $1/2$.]{
    	\psfig{file=bg7fig3c.eps,height=1.5in}}
     \hspace{0.3 in}
  \subfigure[$C(\om)$ for $\om$ close to $1/2$.]{
    	\psfig{file=bg7fig3d.eps,height=1.5in}}
  \subfigure[$C_{\rho}(\om)$ for $\om$ close to $1/3$.]{
    	\psfig{file=bg7fig3e.eps,height=1.5in}}
     \hspace{0.3 in}
  \subfigure[$C(\om)$ for $\om$ close to $1/3$.]{
    	\psfig{file=bg7fig3f.eps,height=1.5in}}
  	\caption{Plot of $C_{\rho}(\om)=\log\rho(\om)+2B(\om)$ and
	$C(\om)=\log\epsc(\om)+B(\om)$ versus $\om$
	for the sequence $\om_{n_{k}}$ listed in tables
	\ref{tab:1} (plots a and b), \ref{tab:3} (plots c and d)
	and \ref{tab:4} (plots e and f).
        }
        \label{fig:3}
\end{figure}

If we try to plot $C_{\rho}(\om)$ for the values of $\om$
close to $0$ listed in table \ref{tab:1}, and use the values of
$\rho(\om)$ in table \ref{tab:11} we find the behavior represented in
figure \ref{fig:3}a, which also support the smoothness conjecture.
Analogously, if we plot $C(\om)$ for the same set of
values of $\om$, by using the values of $\eps(\om)$ listed
in table \ref{tab:1}, we find the behavior in figure \ref{fig:3}b.

In figures \ref{fig:3}c and \ref{fig:3}d we represent the
functions $C_{\rho}(\om)$ and $C(\om)$ for the values
of $\om$ close to $1/2$ listed in table \ref{tab:3}, and
in figures \ref{fig:3}e and \ref{fig:3}f we represent the
functions $C_{\rho}(\om)$ and $C(\om)$ for the values
of $\om$ close to $1/3$ listed in table \ref{tab:4}.

While the variation in the case of $\epsc(\om)$ is larger,
all plots support the conjecture of a function which is
not only bounded but also H\"older continuous close to the resonances,
but of course a deeper numerical investigation is needed in order
to draw more quantitative deductions.

To conclude, we note also that the value $\beta=2$, which holds for the
radius of convergence, is found both in \cite{CL} and in \cite{CJ}
for the critical function, if the computations are made
without requiring a high precision (this is attributed to
the low precision in \cite{CL} and to the truncations and
approximations due to the numerical implementation of the
renormalization group method in \cite{CJ}). We find this phenomenon at
least very curious: it would be in fact quite interesting to understand
why truncation and approximation errors in the numerical computations
give a different value of $\beta$ (and exactly the one holding for the
radius of convergence) instead than just gibberish.

\section{Conclusions}\label{sect:conclusions}

We conclude by some general remarks about the
advances which have been made and the conclusions which can be
drawn from our analysis. 

\begin{enumerate}

\item The numerical results of \cite{CL} have been improved by an order
of magnitude, both in the size of $\Delta\om$ and in the order of the
numerical errors. Moreover, an heuristic argument providing corrections
to the leading order has been given: the analysis of the numerical
data, taking into account the conjectured form of the corrections,
supports both \eqref{eq:interpol2} with $\beta=1$ and the continuity of
the function $C(\om)$. A stronger support would require getting closer
to the resonances, and considering more resonances rather than just
$0/1$, $1/2$ and $1/3$: all these actions are clearly feasible but
would require significantly more computer time, which has already
reached the order of several CPU years on Compaq Alpha computers for the
calculations of the present paper. 

\item While a reasonably complete analysis of what happens for
sequences of rotation numbers which are \emph{not} the best ones cannot
be practically done, the study of a sequence like $[n,20,1^{\infty}]$
shows that the behavior  of the critical function along this sequence
is the same that along the ``best'' sequence $[n,1^{\infty}]$. It is
likely that this holds, in the limit of long periodic orbits, at least
for all sequences of noble numbers tending to a rational value. 
We could investigate very few non-noble sequences and no sequence
at all made by something different than quadratic irrationals (which
have measure $0$). For a sequence like $[n,2^{\infty}]$ the values of 
$\epsc(\om)$ seem to be comparable to the ones of the ``best'' sequence
quoted above. 

\item The study of the functions $C(\om)=\log\epsc(\om)+
B(\om)$ and $C_{\rho}(\om)=\log\rho(\om)+2B(\om)$ seem
to suggest that they depend smoothly on $\om$.
In general such functions look as continuous in their
arguments, as also the comparison between
the two sequences $[n,1^{\infty}]$ and $[20,n,1^{\infty}]$ seem
to support. This is in contrast with the conclusions made
in \cite{LFLM}, where doubts were raised about
continuity of the function $\log\epsc(\om)+\beta B(\om)$
for any choice of $\beta$.

\item Some interesting conclusions can be drawn for the behavior of
the critical residues. It appears that a limit value is obtained only
when $\om$ is of constant type,  and this limit seems to depend only on
the ``tail'' of the expansion. 
If $\om$ is not of that form but still a quadratic irrational,
then the  sequence of the partial quotients is eventually periodic. In
this case it appears  that $\RR_{k}(\epsc(\om))$, for large $k$,
approaches a periodic sequence of values with the same period of the
partial quotients. 
If $\om$ is not a quadratic irrational, then it is
difficult to draw any conclusion at all. If the partial quotients are
bounded, then the critical residues seem to be bounded away from zero
and oscillating in an apparent random way, but bounded; for unbounded
sequences of partial quotients no numerical data at all could be
obtained. We note that such a scenario is consistent with that
arising within the renormalization group approach as described,
for instance, in \cite{CJ}.

\end{enumerate}

\section*{Acknowledgments}\label{sect:ack}

We thank Silio D'Angelo (INFN sez. Roma 2) and the
Department of Physics of the University of Rome ``La Sapienza''
for providing us some of the computing resources.
All calculations have been done on Compaq Alpha computers using 
Fortran 90 and Mathematica.


\section*{Tables}

\begin{table}[H]
\caption{Values of the Bryuno function, of the critical
function and of the running slopes $A_{k}=A(\om_{n_{k}},\om_{n_{k-1}})$
corresponding to a finite sequence of rotation numbers
$\om_{n_{k}}=1/(n_{k}+\ga)=[n_{k},1^{\infty}]$.
The error on $\eps_{c}(\om_{n_{k}})$ is of 1 unit on the last digit,
and the corresponding slopes are computed with consistent accuracy.}
\begin{center}
\setlength\tabcolsep{5pt}
\vrule
\begin{tabular}{llllll}
\hline\noalign{\smallskip}
        $k$ & $\om_{n_{k}}$ & $B(\om_{n_{k}})$ &
        $\epsc(\om_{n_{k}})$ & $A_{k}$\\
\noalign{\smallskip}
\hline
\noalign{\smallskip}
 1   & $  [500,1^{\infty}]$ &  6.21836 & 0.016585    & \\
 2   & $  [700,1^{\infty}]$ &  6.55376 & 0.0121005   & 0.9399$\pm$0.0002 \\
 3   & $ [1000,1^{\infty}]$ &  6.90963 & 0.0086401   & 0.9465$\pm$0.0001 \\
 4   & $ [2000,1^{\infty}]$ &  7.60184 & 0.0044599   & 0.9553$\pm$0.0001\\
 5   & $ [4000,1^{\infty}]$ &  8.29452 & 0.0022854   & 0.9652$\pm$0.0001 \\
 6   & $ [7000,1^{\infty}]$ &  8.85393 & 0.0013265   & 0.9724$\pm$0.0002 \\
 7   & $[10000,1^{\infty}]$ &  9.21053 & 0.00093627  & 0.9770$\pm$0.0002 \\
 8   & $[12000,1^{\infty}]$ &  9.39284 & 0.00078320  & 0.9793$\pm$0.0001 \\
 9   & $[15000,1^{\infty}]$ &  9.61593 & 0.00062927  & 0.9808$\pm$0.0001 \\
10   & $[18000,1^{\infty}]$ &  9.79823 & 0.00052610  & 0.9823$\pm$0.0002 \\
11   & $[20000,1^{\infty}]$ &  9.90358 & 0.00047433  & 0.9833$\pm$0.0004 \\
12   & $[25000,1^{\infty}]$ & 10.12671 & 0.00038081  & 0.9842$\pm$0.0002 \\
13   & $[30000,1^{\infty}]$ & 10.30902 & 0.00031816  & 0.9859$\pm$0.0003 \\
14   & $[40000,1^{\infty}]$ & 10.59668 & 0.00023955  & 0.9865$\pm$0.0003 \\
15   & $[50000,1^{\infty}]$ & 10.81982 & 0.000192161 & 0.9879$\pm$0.0003 \\
16   & $[60000,1^{\infty}]$ & 11.00213 & 0.000160443 & 0.9895$\pm$0.0002 \\
\hline
\end{tabular}
\hspace{-0.1cm}\vrule
\end{center}
\label{tab:1}
\end{table}

\begin{table}[H]
\caption{Values of the Bryuno function, of the critical
function and of the running slopes $A_{k}=A(\om_{n_{k}},\om_{n_{k-1}})$
corresponding to a finite sequence of rotation numbers
$\om_{n_{k}}=1/(n_{k}+1/(20 +\ga))=[n_{k},20,1^{\infty}]$.
The error on $\eps_{c}(\om_{n_{k}})$ is of 1 unit on the last digit,
and the corresponding slopes are computed with consistent accuracy.}
\begin{center}
\setlength\tabcolsep{5pt}
\vrule
\begin{tabular}{llllll}
\hline\noalign{\smallskip}
        $k$ & $\om_{n_{k}}$ & $B(\om_{n_{k}})$ &
        $\epsc(\om_{n_{k}})$ & $A_{k}$ \\
\noalign{\smallskip}
\hline
\noalign{\smallskip}
 1 & $ [500,20,1^{\infty}]$ &  6.22088 & 0.016303   &\\
 2 & $ [700,20,1^{\infty}]$ &  6.55556 & 0.011926   & 0.9341$\pm$0.0004 \\
 3 & $[1000,20,1^{\infty}]$ &  6.91089 & 0.008535   & 0.9415$\pm$0.0006 \\
 4 & $[2000,20,1^{\infty}]$ &  7.60247 & 0.004421   & 0.9512$\pm$0.0005 \\
 5 & $[4000,20,1^{\infty}]$ &  8.29483 & 0.002271   & 0.962 $\pm$0.001  \\
\hline
\end{tabular}
\hspace{-0.1cm}\vrule
\end{center}
\label{tab:2}
\end{table}

\begin{table}[H]
\centering
\caption{Values of the Bryuno function, of the critical
function and of the running slopes $A_{k}=A(\om_{n_{k}},\om_{n_{k-1}})$
corresponding to a finite sequence of rotation numbers
$\om_{n_{k}}=1/(2+1/(n_{k}+\ga))=[2,n_{k},1^{\infty}]$.
The error on $\eps_{c}(\om_{n_{k}})$ is of 1 unit on the last digit,
and the corresponding quantities are computed with consistent accuracy.}
\begin{center}
\setlength\tabcolsep{5pt}
\vrule
\begin{tabular}{lllllll}
\hline\noalign{\smallskip}
        $k$ & $\om_{n_{k}}$ & $B(\om_{n_{k}})$ &
        $\epsc(\om_{n_{k}})$ & $A_{k}$ \\
\noalign{\smallskip}
\hline
\noalign{\smallskip}
 1   &   $[2,500,1^{\infty}]$ & 3.80022 & 0.12872   &\\
 2   &   $[2,700,1^{\infty}]$ & 3.96840 & 0.109967  & 0.9362$\pm$0.0005\\
 3   &  $[2,1000,1^{\infty}]$ & 4.14674 & 0.092932  & 0.9438$\pm$0.0001\\
 4   &  $[2,2000,1^{\infty}]$ & 4.49337 & 0.066777  & 0.9535$\pm$0.0001\\
 5   &  $[2,4000,1^{\infty}]$ & 4.84001 & 0.047805  & 0.9642$\pm$0.0001\\
 6   &  $[2,7000,1^{\infty}]$ & 5.11987 & 0.036420  & 0.9720$\pm$0.0002\\
 7   & $[2,10000,1^{\infty}]$ & 5.29823 & 0.030598  & 0.9766$\pm$0.0003\\
 8   & $[2,13000,1^{\infty}]$ & 5.42943 & 0.026909  & 0.9792$\pm$0.0006\\
 9   & $[2,17000,1^{\infty}]$ & 5.56357 & 0.023591  & 0.9810$\pm$0.0006\\
10   & $[2,20000,1^{\infty}]$ & 5.64484 & 0.021780  & 0.983 $\pm$0.001 \\
\hline
\end{tabular}
\hspace{-0.1cm}\vrule
\end{center}
\label{tab:3}
\end{table}

\begin{table}[H]
\centering
\caption{Values of the Bryuno function, of the critical
function and of the running slopes $A_{k}=A(\om_{n_{k}},\om_{n_{k-1}})$ 
corresponding to a finite sequence of rotation numbers
$\om_{n_{k}}=1/(3+1/(n_{k}+\ga))=[3,n_{k},1^{\infty}]$.
The error on $\eps_{c}(\om_{n_{k}})$ is of 1 unit on the last digit,
and the corresponding quantities are computed with consistent accuracy.}
\begin{center}
\setlength\tabcolsep{5pt}
\vrule
\begin{tabular}{llllll}
\hline\noalign{\smallskip}
        $k$ & $\om_{n_{k}}$ & $B(\om_{n_{k}})$ &
        $\epsc(\om_{n_{k}})$ & $A_{k}$\\
\noalign{\smallskip}
\hline
\noalign{\smallskip}
 1   &   $[3,500,1^{\infty}]$ & 3.17069 & 0.244787  & \\
 2   &   $[3,700,1^{\infty}]$ & 3.28264 & 0.22044   & 0.9358$\pm$0.0001\\
 3   &  $[3,1000,1^{\infty}]$ & 3.40139 & 0.197080  & 0.9433$\pm$0.0001\\
 4   &  $[3,2000,1^{\infty}]$ & 3.63230 & 0.158153  & 0.9529$\pm$0.0001\\
 5   &  $[3,4000,1^{\infty}]$ & 3.86330 & 0.126588  & 0.9637$\pm$0.0001\\
 6   &  $[3,7000,1^{\infty}]$ & 4.04983 & 0.105608  & 0.9715$\pm$0.0001\\
 7   & $[3,10000,1^{\infty}]$ & 4.16872 & 0.094035  & 0.9763$\pm$0.0002\\
 8   & $[3,13000,1^{\infty}]$ & 4.25617 & 0.086319  & 0.9787$\pm$0.0003\\
 9   & $[3,17000,1^{\infty}]$ & 4.34559 & 0.079072  & 0.9807$\pm$0.0003\\
10   & $[3,20000,1^{\infty}]$ & 4.39977 & 0.074973  & 0.9826$\pm$0.0005\\
\hline
\end{tabular}
\hspace{-0.1cm}\vrule
\end{center}
\label{tab:4}
\end{table}

\begin{table}[H]
\centering
\caption{Critical residues $\RR_{\infty}(\om)$
for some rotation numbers $\om$.
The error on $\eps_{c}(\om_{n_{k}})$ is of 1 unit on the last digit.}
\begin{center}
\setlength\tabcolsep{5pt}
\vrule
\begin{tabular}{llll}
\hline\noalign{\smallskip}
        $\om$ & $\eps_{c}(\om)$ & $\RR_{\infty}(\om)$ \\
\noalign{\smallskip}
\hline
\noalign{\smallskip}
         $[1^{\infty}]$     &  0.971635406 &  0.250088   \\
         $[2^{\infty}]$     &  0.957445408 &  0.2275138  \\
         $[3^{\infty}]$     &  0.890863502 &  0.202230   \\
         $[4^{\infty}]$     &  0.80472544  &  0.17923    \\
         $[10,2^{\infty}]$  &  0.481985986 &  0.22751    \\
         $[1,3,2^{\infty}]$ &  0.829500533 &  0.22751    \\
         $[7,3^{\infty}]$   &  0.615071885 &  0.2022     \\
         $[1,2,4^{\infty}]$ &  0.86423037  &  0.1792     \\
\hline
\end{tabular}
\hspace{-0.1cm}\vrule
\end{center}
\label{tab:5}
\end{table}

\begin{table}[H]
\centering
\caption{Residues of critical periodic orbits for 
$\om=\sqrt{3}-1 = [1,2,1,2,1,2,\dots]$.}
\begin{center}
\setlength\tabcolsep{5pt}
\vrule
\begin{tabular}{lll}
\hline\noalign{\smallskip}
        $\eps_{c}(\om)$ & $0.876067426$ \\
\noalign{\smallskip}
\hline
\noalign{\smallskip}
        approximant          &  residue \\
\noalign{\smallskip}
\hline
\noalign{\smallskip}
3/4              & 0.24871   \\
8/11             & 0.18612   \\
11/15            & 0.25216   \\
30/41            & 0.18516   \\
41/56            & 0.25275   \\
112/153          & 0.18493   \\
153/209          & 0.25288   \\
418/571          & 0.18487   \\
\hline
\end{tabular}
\hspace{-0.1cm}\vrule
\hspace{1cm}
\vrule
\begin{tabular}{lll}
\hline\noalign{\smallskip}
        $\eps_{c}(\om)$ & $0.876067426$ \\
\noalign{\smallskip}
\hline
\noalign{\smallskip}
        approximant          &  residue \\
\noalign{\smallskip}
\hline
\noalign{\smallskip}
571/780          & 0.25291   \\
1560/2131        & 0.18486   \\
2131/2911        & 0.25292   \\
5822/7953        & 0.18485   \\
7953/10864       & 0.25292   \\
21728/29681      & 0.18485   \\
29681/40545      & 0.25292   \\
81090/110771     & 0.18486   \\
\hline
\end{tabular}
\hspace{-0.1cm}\vrule
\end{center}
\label{tab:6}
\end{table}

\begin{table}[H]
\centering
\caption{Residues of critical periodic orbits for 
$\om=(\sqrt{3}-1)/2 = [2,1,2,1,2,1,\dots]$.}
\begin{center}
\setlength\tabcolsep{5pt}
\vrule
\begin{tabular}{lll}
\hline\noalign{\smallskip}
        $\eps_{c}(\om)$ & $0.9402827$ \\
\noalign{\smallskip}
\hline
\noalign{\smallskip}
        approximant          &  residue \\
\noalign{\smallskip}
\hline
\noalign{\smallskip}
3/8              &  0.19574  \\
4/11             &  0.24746  \\
11/30            &  0.18763  \\
15/41            &  0.25145  \\
41/112           &  0.18556  \\
56/153           &  0.25254  \\
153/418          &  0.18503  \\
209/571          &  0.25282  \\
\hline
\end{tabular}
\hspace{-0.1cm}\vrule
\hspace{1cm}
\vrule
\begin{tabular}{lll}
\hline\noalign{\smallskip}
        $\eps_{c}(\om)$ & $0.9402827$ \\
\noalign{\smallskip}
\hline
\noalign{\smallskip}
        approximant          &  residue \\
\noalign{\smallskip}
\hline
\noalign{\smallskip}
571/1560         &  0.18490  \\
780/2131         &  0.25290  \\
2131/5822        &  0.18486  \\
2911/7953        &  0.25292  \\
7953/21728       &  0.18486  \\
10864/29681      &  0.25292  \\
29681/81090      &  0.18486  \\
40545/110771     &  0.25293  \\
\hline
\end{tabular}
\hspace{-0.1cm}\vrule
\end{center}
\label{tab:7}
\end{table}

\begin{table}[H]
\centering
\caption{Residues of critical periodic orbits for 
$\om=\sqrt{5/2}-1 = [1,1,2,1,1,2,1,1,2\dots]$.}
\begin{center}
\setlength\tabcolsep{5pt}
\vrule
\begin{tabular}{lll}
\hline\noalign{\smallskip}
        $\eps_{c}(\om)$ & $0.9402827$ \\
\noalign{\smallskip}
\hline
\noalign{\smallskip}
        approximant          &  residue \\
\noalign{\smallskip}
\hline
\noalign{\smallskip}
3/5                  & 0.2242  \\
4/7                  & 0.2639  \\
7/12                 & 0.2278  \\
18/31                & 0.2222  \\
25/43                & 0.2660  \\
43/74                & 0.2270  \\
111/191              & 0.2227  \\
154/265              & 0.2656  \\
265/456              & 0.2272  \\
\hline
\end{tabular}
\hspace{-0.1cm}\vrule
\hspace{1cm}
\vrule
\begin{tabular}{lll}
\hline\noalign{\smallskip}
        $\eps_{c}(\om)$ & $0.9402827$ \\
\noalign{\smallskip}
\hline
\noalign{\smallskip}
        approximant          &  residue \\
\noalign{\smallskip}
\hline
\noalign{\smallskip}
684/1177             & 0.2226  \\
949/1633             & 0.2656  \\
1633/2810            & 0.2271  \\
4215/7253            & 0.2227  \\
5848/10063           & 0.2656  \\
10063/17316          & 0.2271  \\
25974/44695          & 0.2227  \\
36037/62011          & 0.2656  \\
62011/106706         & 0.2272  \\
\hline
\end{tabular}
\hspace{-0.1cm}\vrule
\end{center}
\label{tab:8}
\end{table}

\begin{table}[H]
\centering
\caption{Radius of convergence for some values of the
rotation number $\om$ close to $1/2$ and slopes 
$A_{k}'=A'(\om_{n_{k}},\om_{n_{k-1}})$. The value $\rho_{1}(\om)$
is given by the formula \eqref{eq:rhoapprox}, while $\rho_{P}(\om)$
is the value obtained numerically by using Pad\'{e} approximants.
the two values for the slopes correspond to the
values $\rho_{1}(\om)$ and $\rho_{P}(\om)$, respectively.
One has $\eta=|\om-1/2|$.}
\begin{center}
\setlength\tabcolsep{5pt}
\vrule
\begin{tabular}{lllllll}
\hline\noalign{\smallskip}
	$k$ & $\om_{n_{k}}$ & $\eta$
	& $\rho_{1}(\om_{n_{k}})$& $\rho_{P}(\om_{n_{k}})$ & 
	$A_{k}'$\\
\noalign{\smallskip}
\hline
\noalign{\smallskip}
1 & $[2,10,1^{\infty}]$ & 0.0224860 & 0.51409 & 0.51052 & \\
2 & $[2,12,1^{\infty}]$ & 0.0190577 & 0.43571 & 0.43355 & 2.19667/2.17013 \\
3 & $[2,15,1^{\infty}]$ & 0.0155106 & 0.35462 & 0.35352 & 2.14426/2.12464 \\
4 & $[2,20,1^{\infty}]$ & 0.0118382 & 0.27066 & 0.27024 & 2.09658/2.08449 \\
5 & $[2,30,1^{\infty}]$ & 0.0080339 & 0.18368 & 0.18361 & 2.05439/2.04822 \\
6 & $[2,40,1^{\infty}]$ & 0.0060801 & 0.13901 & 0.13902 & 2.02821/2.02484 \\
6 & $[2,50,1^{\infty}]$ & 0.0048906 & 0.11181 & 0.11184 & 2.01612/2.01480 \\
\hline
\end{tabular}
\hspace{-0.1cm}\vrule
\end{center}
\label{tab:9}
\end{table}

\begin{table}[H]
\centering
\caption{Radius of convergence for some values of the
rotation number $\om$ close to $1/3$ and slopes 
$A_{k}'=A'(\om_{n_{k}},\om_{n_{k-1}})$. The value $\rho_{1}(\om)$
is given by the formula \eqref{eq:rhoapprox}, while $\rho_{P}(\om)$
is the value obtained numerically by using Pad\'{e} approximants;
the two values for the slopes correspond to the
values $\rho_{1}(\om)$ and $\rho_{P}(\om)$, respectively.
One has $\eta=|\om-1/3|$.}
\begin{center}
\setlength\tabcolsep{5pt}
\vrule
\begin{tabular}{lllllll}
\hline\noalign{\smallskip}
	$k$ & $\om_{n_{k}}$ & $\eta$ &
	$\rho_{1}(\om_{n_{k}})$ & $\rho_{P}(\om_{n_{k}})$ &
	$A_{k}'$ \\
\noalign{\smallskip}
\hline
\noalign{\smallskip}
1 & $[3,10,1^{\infty}]$  & 0.0101459 & 0.62329 & 0.61993 & \\
2 & $[3,12,1^{\infty}]$  & 0.0085791 & 0.55734 & 0.55524 & 2.28295/2.24934 \\
3 & $[3,13,1^{\infty}]$  & 0.0079642 & 0.53038 & 0.52858 & 2.23762/2.22067 \\
4 & $[3,20,1^{\infty}]$  & 0.0053033 & 0.40444 & 0.40400 & 2.17212/2.15360 \\
5 & $[3,30,1^{\infty}]$  & 0.0035899 & 0.31180 & 0.31182 & 2.09982/2.09051 \\
6 & $[3,40,1^{\infty}]$  & 0.0027132 & 0.25871 & 0.25872 & 2.06311/2.06339 \\
7 & $[3,50,1^{\infty}]$  & 0.0021807 & 0.22364 & 0.22360 & 2.04490/2.04795 \\
8 & $[3,100,1^{\infty}]$ & 0.0011006 & 0.14177 & 0.14179 & 2.02455/2.02313 \\
9 & $[3,200,1^{\infty}]$ & 0.0005529 & 0.08959 & 0.08961 & 2.00902/2.00866 \\
\hline
\end{tabular}
\hspace{-0.1cm}\vrule
\end{center}
\label{tab:10}
\end{table}

\begin{table}[H]
\centering
\caption{Values of the radius of
convergence and of the slopes $A_{k}'=A'(\om_{n_{k}},\om_{n_{k-1}})$
corresponding to a finite sequence of rotation numbers
$\om_{n_{k}}=1/(n_{k}+\ga)=[n_{k},1^{\infty}]$.
The radius of convergence is computed
with the formula \eqref{eq:rhoapprox}.}
\begin{center}
\setlength\tabcolsep{5pt}
\vrule
\begin{tabular}{llll}
\hline\noalign{\smallskip}
        $k$ & $\om_{n_{k}}$ &
        $\rho(\om_{n_{k}})$ &
	$A_{k}'$\\
\noalign{\smallskip}
\hline
\noalign{\smallskip}
 1   &   $[500,1^{\infty}]$ & 0.000130355     & \\
 2   &   $[700,1^{\infty}]$ & 0.0000665545    & 2.0042837 \\
 3   &  $[1000,1^{\infty}]$ & 0.000032629     & 2.0030298 \\
 4   &  $[2000,1^{\infty}]$ & 0.00000816229   & 2.0018183 \\
 5   &  $[4000,1^{\infty}]$ & 0.0000020412    & 2.0009090 \\
 6   &  $[7000,1^{\infty}]$ & 0.000000666603  & 2.0004825 \\
 7   & $[10000,1^{\infty}]$ & 0.000000326653  & 2.0003028 \\
 8   & $[12000,1^{\infty}]$ & 0.000000226847  & 2.0002303 \\
 9   & $[15000,1^{\infty}]$ & 0.000000145185  & 2.0001882 \\
10   & $[18000,1^{\infty}]$ & 0.000000100824  & 2.0001536 \\
11   & $[20000,1^{\infty}]$ & 0.0000000816683 & 2.0001329 \\
12   & $[25000,1^{\infty}]$ & 0.0000000522684 & 2.0001129 \\
13   & $[30000,1^{\infty}]$ & 0.0000000362978 & 2.0000921 \\
14   & $[40000,1^{\infty}]$ & 0.0000000204177 & 2.0000730 \\
15   & $[50000,1^{\infty}]$ & 0.0000000130674 & 2.0000565 \\
\hline
\end{tabular}
\hspace{-0.1cm}\vrule
\end{center}
\label{tab:11}
\end{table}

\begin{table}[H]
\centering
\caption{Values of the radius of
convergence and of the slopes $A_{k}'=A'(\om_{n_{k}},\om_{n_{k-1}})$
corresponding to a finite sequence of rotation numbers
$\om_{n_{k}}=1/(2+1/(n_{k}+\ga))=[2,n_{k},1^{\infty}]$.
The radius of convergence is computed
with the formula \eqref{eq:rhoapprox}.}
\begin{center}
\setlength\tabcolsep{5pt}
\vrule
\begin{tabular}{llll}
\hline\noalign{\smallskip}
        $k$ & $\om_{n_{k}}$ &
        $\rho(\om_{n_{k}})$ &
	$A_{k}'$ \\
\noalign{\smallskip}
\hline
\noalign{\smallskip}
 1   &   $[2,500,1^{\infty}]$ & 0.011405915 & \\
 2   &   $[2,700,1^{\infty}]$ & 0.008152279 & 1.9968638 \\
 3   &  $[2,1000,1^{\infty}]$ & 0.005709327 & 1.9973651 \\
 4   &  $[2,2000,1^{\infty}]$ & 0.002856258 & 1.9980597 \\
 5   &  $[2,4000,1^{\infty}]$ & 0.001428528 & 1.9987793 \\
 6   &  $[2,7000,1^{\infty}]$ & 0.000816400 & 1.9992292 \\
 7   & $[2,10000,1^{\infty}]$ & 0.000571507 & 1.9994593 \\
 8   & $[2,13000,1^{\infty}]$ & 0.000439632 & 1.9995765 \\
 9   & $[2,17000,1^{\infty}]$ & 0.000336196 & 1.9996572 \\
10   & $[2,20000,1^{\infty}]$ & 0.000285770 & 1.9997123 \\
\hline
\end{tabular}
\hspace{-0.1cm}\vrule
\end{center}
\label{tab:12}
\end{table}

\begin{table}[H]
\centering
\caption{Values of the radius of
convergence and of the slopes $A_{k}'=A'(\om_{n_{k}},\om_{n_{k-1}})$
corresponding to a finite sequence of rotation numbers
$\om_{n_{k}}=1/(3+1/(n_{k}+\ga))=[3,n_{k},1^{\infty}]$.
The error on $\eps_{c}(\om_{n_{k}})$ is of 1 unit on the last digit,
and the corresponding logarithm is computed with consistent accuracy.}
\begin{center}
\setlength\tabcolsep{5pt}
\vrule
\begin{tabular}{lllllll}
\hline\noalign{\smallskip}
        $k$ & $\om_{n_{k}}$ &
        $\rho(\om_{n_{k}})$ & 
	$A_{k}'$ \\
\noalign{\smallskip}
\hline
\noalign{\smallskip}
 1   &   $[3,500,1^{\infty}]$ & 0.04873028  & \\
 2   &   $[3,700,1^{\infty}]$ & 0.03895268  & 2.0004598 \\
 3   &  $[3,1000,1^{\infty}]$ & 0.03071760  & 2.0000489 \\
 4   &  $[3,2000,1^{\infty}]$ & 0.01935701  & 1.9997910 \\
 5   &  $[3,4000,1^{\infty}]$ & 0.01219611  & 1.9997289 \\
 6   &  $[3,7000,1^{\infty}]$ & 0.00839894  & 1.9997744 \\
 7   & $[3,10000,1^{\infty}]$ & 0.00662168  & 1.9998204 \\
 8   & $[3,13000,1^{\infty}]$ & 0.00555920  & 1.9998501 \\
 9   & $[3,17000,1^{\infty}]$ & 0.00464887  & 1.9998731 \\
10   & $[3,20000,1^{\infty}]$ & 0.00417153  & 1.9998900 \\
\hline
\end{tabular}
\hspace{-0.1cm}\vrule
\end{center}
\label{tab:13}
\end{table}

\end{document}